\documentclass[11pt]{amsart}
\usepackage{amsmath,amsthm,epsf,amscd,amssymb,verbatim}

\newtheorem{thm}{Theorem}[section]
\newtheorem{lem}[thm]{Lemma}

\newtheorem{prop}[thm]{Proposition}

\theoremstyle{definition}
\newtheorem{dfn}[thm]{Definition}

\theoremstyle{remark}
\newtheorem{remark}[thm]{Remark}

\numberwithin{equation}{section}

\newcommand{\bbr}{\begin{remark}}        
\newcommand{\eer}{\end{remark}}

\font\bbb=msbm10 scaled 1100

\newcommand{\be}{\begin{equation}}
\newcommand{\ee}{\end{equation}}
\newcommand{\bea}{\begin{eqnarray}}
\newcommand{\eea}{\end{eqnarray}}
\newcommand{\bmini}{\footnotesize\begin{center}\begin{minipage}{5.5in}}
\newcommand{\emini}{\end{minipage}\end{center}\normalsize}
\newcommand{\pf}{{\em Proof: }}

\newcommand{\real}{\mbox{\bbb R}}

\newcommand{\eg}{{\em e.g.}}
\newcommand{\ie}{{\em i.e.}}

\newcommand{\mathspace}{\;\;\;\;\;}
\newcommand{\del}{\partial}

\newcommand{\dbyd}[2]{\frac{\displaystyle\partial #1}{\displaystyle\partial #2}}

\newcommand{\norm}[1]{\left\|{#1}\right\|}
\newcommand{\grad}{\nabla}

\newcommand{\curl}{\nabla\times}

\newcommand{\Lie}{{\mathcal L}}
\newcommand{\Index}{{\mbox{\bbb I}}}
\newcommand{\Ind}{{\mbox{\rm Ind}}}
\newcommand{\sign}{{\mbox{\rm Sign}}}
\begin{document}

\title{AN INDEX FOR CLOSED ORBITS IN BELTRAMI FIELDS}
\author{John Etnyre}
\address{Department of Mathematics, Stanford University, 
	Palo Alto, CA, 94305}
\thanks{JE supported in part by NSF Grant \# DMS-9705949.}
\email{{\tt etnyre@math.stanford.edu}}

\author{Robert Ghrist}
\address{School of Mathematics, Georgia Institute of Technology, Atlanta, 
	GA 30332}
\thanks{RG supported in part by NSF Grant \# DMS-9971629. The writing 
of this paper was facilitated by the hospitality of the Isaac Newton 
Institute at Cambridge University.}
\email{{\tt ghrist@math.gatech.edu}}


\begin{abstract}
We consider the class of Beltrami fields (eigenfields of the 
curl operator) on three-dimensional Riemannian solid tori: such 
vector fields arise as steady incompressible inviscid fluids 
and plasmas. Using techniques from contact geometry, we construct 
an integer-valued index for detecting closed orbits
in the flow which are topologically inessential (they have
winding number zero with respect to the solid torus). This index
is independent of the Riemannian structure, and is 
computable entirely from a $C^1$ approximation to the 
vector field on any meridional disc of the solid torus. 
\end{abstract}

\maketitle

\section{Introduction and summary}

Consider the class of {\em Beltrami fields} --- the 
volume-preserving eigenfields of the curl operator. 
Such vector fields are the source of numerous interesting 
phenomena in inviscid fluids and plasmas. For example, 
Beltrami fields are the only steady three-dimensional Euler flows 
which admit chaotic Lagrangian dynamics. Beltrami fields are
also common approximations to the magnetic field lines in 
large-scale structures within the solar corona. 

Despite their importance and inherent intricacy, very little 
is known about the dynamics of Beltrami fields apart from 
numerical simulation \cite{Hen66,Dom+86}
and Melnikov analyses of near-integrable Beltrami fields
\cite{Gau85,ZKBH93,GT95,Chi96,HZD98} -- an important but 
extremely small class of solutions. 
We consider the subtle problem of understanding how much
and what kinds of dynamics Beltrami fields are forced to
possess given the underlying topological features of the 
fluid domain. 

In a series of papers 
\cite{EG:I,EG:II,EG:III}, the authors develop techniques
for determining forced behaviors in steady inviscid fluids 
via the topology of {\em contact structures}, 
the odd-dimensional analogues of symplectic structures
(see, \eg, \cite{Aeb94,ET97} for an introduction to contact geometry). 
In this paper, we give an application of these 
techniques to Beltrami fields on solid tori. One of the 
features of our topological approach is that it is independent 
of the Riemannian metric and is furthermore robust 
with respect to perturbations of the vector field,
without resorting to any hyperbolicity or nondegeneracy 
assumptions usually required to preserve closed orbits.

We restrict attention to Beltrami fields on Riemannian solid 
tori, such as would occur in the case of a force-free plasma in a 
containment device. In \cite{EG:II}, it is shown using techniques 
from contact topology and pseudo-holomorphic curves that 
steady Euler fields on an invariant solid torus {\em always}
possess a closed orbit, independent of the Riemannian
structure and volume form:

\noindent
{\bf Theorem:} (\cite{EG:II}) 
{\em 
Any steady real-analytic solution
to the Euler equations (\ref{eq_euler}) on any invariant 
Riemannian solid torus possesses a closed orbit.
}

Since every Beltrami field is a steady Euler field, the 
result holds true for all Beltrami fields. The restrictive
smoothness assumption is necessary for using singularity
theory arguments for the integrable Euler fields --- in 
the setting of pure Beltrami fields, the techniques are
valid up to smoothness class $C^2$ \cite{Hof93}.

In this paper, we define an integer-valued index for
detecting the presence of {\em contractible} closed 
orbits --- those closed orbits which can be shrunk to a 
point within the solid torus. Thus, we are not primarily
concerned with the class of Beltrami fields which 
possess a cross-section (and hence trivially have a 
closed non-contractible orbit by the Brouwer fixed point
theorem). In cases where no section exists, it is very 
difficult to determine the existence of closed orbits. Indeed, 
the recent examples of fixed-point-free vector fields on a solid 
torus without any periodic orbits (constructible via \cite{KK94}
in the real-analytic case, and by \cite{Kup96} in the $C^1$
volume-preserving case), demonstrate the delicacy of the problem. 

A related scenario to which our results apply is that of a 
Beltrami field on a long tube $D^2\times\real$ which is 
periodic in the third variable. The problem of finding 
contractible closed orbits in the solid torus obtained 
by quotienting out the periodicity is precisely the 
problem of finding an orbit which is closed in the long tube. 

The index we construct is a type of linking number with respect 
to a contact structure --- the so-called {\em self-linking number} 
of a transverse knot, well-known to contact topologists. The data 
required to compute the index is minimal: one needs information 
about the Beltrami field along some (arbitrary) 
meridional disc in the solid torus. The principal contribution of 
this note is to retool the contact-topological index in a form 
which requires no knowledge of the contact structure per se. 
Since contact structures are very stable with respect to perturbation,
the vector field need only be known approximately ($C^1$)
along the disc. Hence, this index can be computed 
numerically with full rigor. 

Sections \ref{sec_beltrami} through \ref{sec_reeb} 
assemble the relevant ingredients. Section
\ref{sec_ot} presents the technical result on contact 
structures used to define 
the index and prove its major properties. The section 
\S\ref{sec_compute} gives a very simple method for
computing this index from a minimal amount of data: one 
need simply know what the vector field $X$ approximately looks like near 
some finite number of points on a meridional disc of the solid torus. 
It is our hope that this computability may allow for utilization of
this index in the analysis of experimental data. 

\section{Beltrami fields on Riemannian manifolds}
\label{sec_beltrami}

Let $M$ be an arbitrary 3-manifold with Riemannian metric $g$ and 
(arbitrary) volume form $\mu$. Given a vector field $X$ on $M$, 
one can consider the dual 1-form $g(X,\cdot)$ to $X$ that pairs 
with a vector $Y$ to give the inner product $g(X,Y)$. In this 
general setting, the {\em curl} of a vector field $X$ on $M$ is the 
unique vector field $\curl X$ satisfying 
\begin{equation}
	\mu({\curl X},\cdot,\cdot) = d(g(X,\cdot)), 
\end{equation}
where $d$ denotes the exterior derivative on forms. 
The curl operator is linear and its $\mu$-preserving 
eigenfields are known as the {\em Beltrami fields}. In 
other words, $X$ is Beltrami if and only if it is volume-preserving 
and $\curl X = \lambda X$ for some constant $\lambda$. 
One can also consider the class of eigenfields with 
scalar fields as eigenvalues: $\curl X = fX$ for $f:M\to\real$.
Our techniques are adaptable to this more general format, 
but we restrict to pure eigenfields here for simplicity.

Beltrami fields arise in several contexts:
\begin{enumerate}
\item
	Beltrami fields are always steady solutions to the 
	Euler equations for an inviscid incompressible fluid
\begin{equation}
\label{eq_euler}
	\dbyd{u}{t} + \nabla_uu = -\grad p \mathspace;\mathspace
	\Lie_u\mu=0 , 
\end{equation}
	where $\nabla_uu$ is the covariant derivative of 
	the velocity field $u$ along itself, and $p:M\to\real$, 
	the pressure function, can be chosen to be $\frac{1}{2}\norm{u}^2$.
	The Lie derivative $\Lie_u\mu$ of the volume form along 
	$u$ vanishing is equivalent to $u$ being divergence-free. 
\item
	Beltrami fields also yield steady solutions to the 
	ideal MHD equations
\[	\dbyd{u}{t} + \nabla_uu = -\grad p + (\curl B)\times B \]

\vspace{-0.15in}
\begin{equation}
	\dbyd{B}{t} - \curl(u\times B) = 0 
\end{equation} 
\vspace{-0.35in}

\[	\Lie_u\mu = \Lie_B\mu = 0 , \]
	where $B$ denotes the magnetic field. In this context, 
	Beltrami fields are known as force-free fields.
\item 
	Beltrami fields are all extrema of the $L^2$ energy 
	functional 
\begin{equation}
	\norm{u}_2 := \frac{1}{2}\int_M\norm{u}^2 d\mu
\end{equation}
	under the action of the volume-preserving diffeomorphism
	group of $M$. Eigenfields of curl 
	having the smallest nonzero eigenvalue globally minimize the 
	energy \cite{Arn86,AK97}.
\end{enumerate}
Beltrami fields also play a role in the analysis of 
the stability of matter \cite{LY86} and in the formation of 
dynamos \cite{CG95}.

The topology and dynamics of Beltrami fields are subtle:
witness the complex dynamics of the well-known {\em ABC 
fields} on the Euclidean 3-torus \cite{Dom+86}. The 
existence of fixed-point-free Beltrami fields on general
Riemannian 3-manifolds is highly nontrivial \cite{EG:I}, 
as is the presence of closed orbits within such
fields \cite{EG:I,EG:II}.

\section{Contact structures and topology}
\label{sec_contact}

Contact structures are the natural complements to Beltrami
fields. Loosely put, a contact structure on a odd-dimensional
manifold $M$ is a hyperplane field which is maximally nonintegrable.
More specifically, on a three-dimensional manifold, a contact
structure is a smoothly-varying plane field (a choice of a 
two-dimensional subspace $\xi_p$ in each tangent space $T_pM$) which cannot be 
stitched together into leaves of a foliation, not even at a point. 
Locally, every contact structure $\xi$ is the kernel of 
a differential 1-form $\alpha$ satisfying the contact condition: 
\begin{equation}
\label{eq_contact}\alpha\wedge d\alpha \neq 0 .
\end{equation}
Otherwise said, $\alpha\wedge d\alpha$ is locally a volume form on $M$. 
Any 1-form satisfying (\ref{eq_contact}) is called a 
contact form. If $\alpha\wedge d\alpha$ is a globally defined 
volume form on $M$, then the contact structure is said to be 
{\em cooriented}: all contact structures which arise in 
connection with Beltrami fields are of necessity cooriented, and 
we will restrict entirely to this category. 

A contact form $\alpha$ on an oriented three-manifold $M$ 
is said to be {\em positive} if the sign of $\alpha\wedge d\alpha$ 
is positive with respect to the orientation of $M$. Otherwise, 
$\alpha$ is said to be {\em negative}. The sign is a property of 
the contact structure and is independent of the defining 1-form. 

Canonical examples of positive contact forms on $\real^3$ include 
$dz+x\,dy$, $dz+r^2\sin r\,d\theta$, and $\cos r\,dz+r\sin r\,d\theta$, 
the latter two being given in cylindrical coordinates. 

Much of the current interest in contact structures arises 
from fairly recent elucidations of their topological 
and dynamical properties (see, \eg, \cite{Aeb94,ET97}). 
Studying contact structures by means of {\em characteristic
foliations} is most fruitful. Given a two-dimensional 
surface $S$ embedded in $M$, the characteristic foliation
of $S$, $S_\xi$, is the [singular] one-dimensional foliation
generated by the intersections of the tangent planes of 
$T_pS$ with the contact planes $\xi_p$ in $T_pM$. For all 
intents and purposes, $S_\xi$ may be thought of as a 
vector field on $S$ generated by $\xi$ (by orienting the 
foliation). The singularities which arise on a characteristic 
foliation are generically saddles or spiral sources/sinks: pure
centers cannot ever appear from a contact structure --- the 
plane field twists too much for this. 

\begin{figure}[ht]
	{\epsfxsize=2.0in\epsfbox{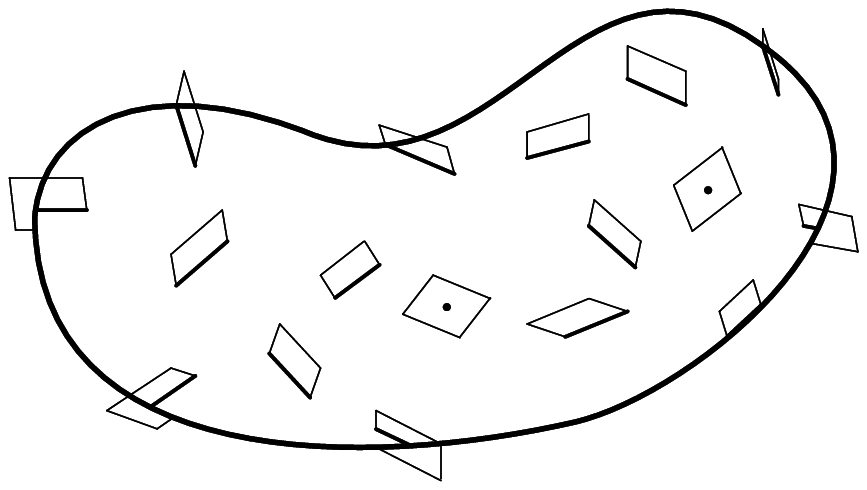}}
	\hspace{0.25in}
	{\epsfxsize=2.0in\epsfbox{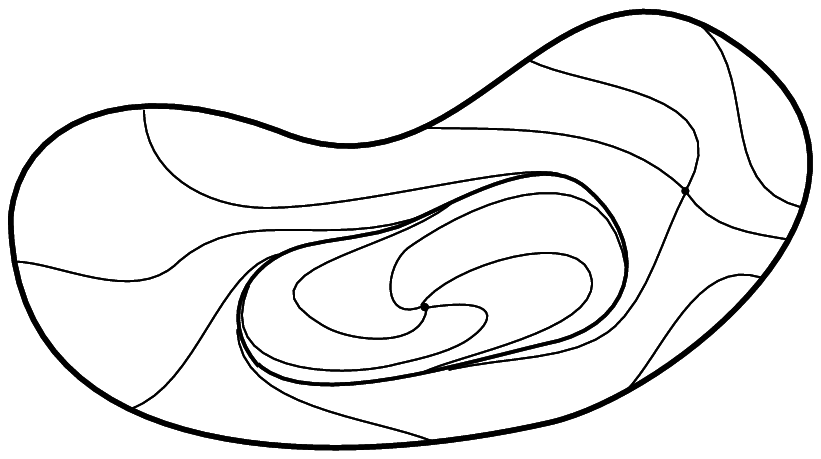}}
	\caption{A contact structure is overtwisted if the induced
	characteristic foliation on some embedded disc possesses 
	a limit cycle.} 
	\label{fig_ot}
\end{figure}

The dynamical properties of $S_\xi$ are 
closely related to the topological classification of 
contact structures. A contact structure $\xi$ is said to 
be {\em overtwisted} if there exists an embedded disc
$D\subset M$ such that $D_\xi$ possesses a limit cycle --- 
a closed orbit along which nearby orbits accumulate [see
Figure~\ref{fig_ot}].  
A contact structure is said to be {\em tight} if there
are no such overtwisted discs anywhere in $M$. The contact
structures for $dz+x\,dy$ and $dz+r^2d\theta$ are tight \cite{Ben83}, 
while that of $\cos r\,dz + r\sin r\,d\theta$ is overtwisted
(\eg, at the disc $\{r \leq 1, z=r^2\}$). 

It is by no means apparent that the above definition is 
at all helpful: but in fact, the entire topological theory
hangs on this dichotomy. Many major questions 
about contact structures are solved in the overtwisted 
category and unknown in the tight category (or, if known, then 
only recently and then by great skill and effort). For example,
while overtwisted contact structures have been completely 
classified up to isotopy \cite{Eli89}, the classification 
of tight structures appears [from what is presently known]
to be delicate at best, intractable at worst \cite{Ben83,Gir94,Kan95, 
Etn97:lens,Hon99a,Hon99b,EH99,EH00}. 

\section{Contact dynamics}
\label{sec_reeb}

The connection between contact structures and Beltrami 
fields is, in the present context, quite straightforward. 
Given any fixed-point-free Beltrami field $X$ on $M$, 
the Beltrami condition states that 
\begin{equation}
	\mu(\lambda X,\cdot,\cdot) = d(g(X,\cdot)) .
\end{equation} 
>From this one can derive the crucial observation that the plane 
field orthogonal to any nonvanishing Beltrami field 
is indeed a contact structure as follows. Witness 
the 1-form $\alpha := g(X,\cdot)$ dual to $X$ via $g$. The kernel of 
this 1-form represents the orthogonal plane field $\xi$ to $X$. This 
form $\alpha$ is a contact form on $M$ since 
\begin{equation}
	\alpha \wedge d\alpha := \lambda g(X,\cdot)\wedge\mu(X,\cdot,\cdot) ,
\end{equation}
which, for $\lambda\neq 0$, is nowhere vanishing, as one can easily 
check by evaluating on local orthogonal coordinate bases of the form 
$(e_1:=X/\norm{X},e_2,e_3)$. 

A little more is in fact true: a Beltrami field annihilates
the exterior derivative of the associated contact form, since
\begin{equation}
	(d\alpha)(X,\cdot) = \mu(X,X,\cdot) = 0 .
\end{equation}
Such vector fields are classical objects known as {\em Reeb 
fields.} The Reeb field of a contact form $\alpha$ is the 
unique vector field $Z$ such that $d\alpha(Z,\cdot)=0$ and 
$\alpha(Z)=1$. We have thus observed that any Beltrami
field $X$ (nonsingular with nonzero eigenvalue) is a Reeb 
field for a contact form, after a possible rescaling to force 
$\alpha(X)=1$. In \cite{EG:I}, a broader version 
of this result was demonstrated: namely, that the class
of nonsingular Beltrami fields on $M$ (up to scaling, 
for any Riemannian structure and volume form) is 
identical to the class of Reeb fields (up to scaling, 
for any contact form). 

This theorem allows one to build ``custom'' solutions to 
the steady Euler equations of very high regularity. For
example, \cite{EG:III} builds a single Beltrami field on a 
Riemannian $\real^3$ which possesses closed orbits of all 
possible knot and link types simultaneously. This can be 
viewed as a rigorous manifestation of the delightful results
of Moffatt on knotting in the Euler equations \cite{Mof85}. 

In the present context, we will use this simple correspondence 
between Beltrami and Reeb fields to import technology from 
contact dynamics. Most specifically, we are interested in 
the utility of the tight/overtwisted dichotomy in describing 
the dynamics of Beltrami fields. One extremely important result in
contact dynamics is the following theorem of Hofer \cite{Hof93}: 
{\em Let $\xi$ be an overtwisted contact 
structure on $M$ a compact 3-manifold without boundary. Then the 
Reeb field of any contact form associated to $\xi$ possesses a
closed orbit, some multiple of which bounds a disc in $M$.}

Loosely speaking, a closed orbit (limit cycle) in the characteristic
foliation of a disc in $M$ implies the existence of a closed 
Reeb orbit which bounds a disc in $M$. The proof relies on the 
delicate techniques of pseudo-holomorphic curves in symplectic 
manifolds. One indication of the implicit nature of the proof 
is that the location of the overtwisted disc has little to no 
correlation with the location of the implicated Reeb orbit. 

It follows from well-known 
properties of pseudo-holomorphic curves that the proof of Hofer's
theorem remains valid for a three-manifold with invariant 
boundary (see \cite{EG:II} for details). 
Thus, for the solid torus, it follows that an 
overtwisted Reeb field must possess a contractible periodic 
orbit (since the fundamental group of the solid torus contains
no elements of finite order except the identity). 

\section{Definition of $\Index$}
\label{sec_ot}

>From Hofer's theorem, then, one way to force a contractible
orbit in a Beltrami field is by finding an overtwisted disc in 
the orthogonal contact structure. This is far from trivial, since
it requires searching for overtwisted discs among all possible
embedded discs in $M$: not a computationally feasible task, 
even if the vector field were known analytically (which, in the 
context of an experimentally generated flow is not generally the 
case).   

The classification of contact structures has been successfully
completed only on a selected class of three-manifolds.
The classification of contact structures on the solid torus 
is quite recent and subtle \cite{Gir00,Hon99a}. In particular, it is 
known that on the solid torus there are tight contact structures which 
are ``stably overtwisted'' --- taking some finite covering space 
of the solid torus and lifting the tight contact structure downstairs
yields an overtwisted contact structure on the cover \cite{Mak98,EG:II}
(see Theorem~\ref{thm_otcover}).
While this is a complication for contact topologists, it is a benefit to 
dynamicists.

\begin{lem}
\label{lem_otcover}
Any Beltrami field transverse to a tight contact structure,
some cover of which is overtwisted, must possess a 
contractible closed orbit. 
\end{lem}
\pf
Assume that $\xi$ is a tight contact structure, some finite 
cover $\tilde\xi$ of which is overtwisted. Then, given any 
Beltrami field $X$ associated to $\xi$, lift this to a Beltrami 
field $\tilde X$ on the overtwisted cover. Applying Hofer's theorem
to the cover implies the existence of a contractible periodic
orbit for $\tilde{X}$; however, since a covering space projection
takes orbits to orbits, the closed orbit upstairs (along with 
the disc that it bounds) {\em must} project to a contractible
closed orbit of the original Beltrami field $X$. 
\qed

Thus, our goal is to effectively determine the existence of 
an overtwisted or stably (with respect to coverings) overtwisted
structure on a solid torus given the least amount of information
about a Beltrami field transverse to it. 

To do so, we recall a common index used in contact topology
(see \cite{Aeb94,Eli93} for an introduction). Given a contact 
structure $\xi$ on a three-manifolds $M$, a simple closed curve (knot) is 
called {\em transverse} if its tangents are everywhere 
transverse to the contact planes. Assume that $\gamma$ is 
an oriented simple closed curve in $M$ which bounds a compact 
oriented surface $\Sigma$ in $M$. Then the {\em self-linking 
number} of $\gamma$ with respect to $\xi$ and $\Sigma$ is 
defined as follows. Choose any vector field $Z$ on a neighborhood 
of $\Sigma$ which has no fixed points and which is always
tangent to $\xi$. That this is possible is a simple 
argument involving the classification of plane bundles. 
Then, flow $\gamma$ for a small amount of time under $Z$ 
to obtain a ``push-off'' curve $\gamma_Z$. The self-linking 
number of $\gamma$, $s\ell k(\gamma)$ is then defined as the 
intersection number of $\gamma_Z$ with $\Sigma$ --- \ie, the 
number of transverse intersections, counted algebraically 
using the orientations. This integer, which can be shown 
to be independent of the vector field $Z$ chosen, is an 
invariant of transverse curves up to isotopy through 
transverse curves \cite{Ben83}. On $S^3$, the self-linking 
number is also independent of the surface $\Sigma$ chosen so long 
as it bounds $\gamma$. 

The following recent result allows for an application of this 
index to Beltrami fields.

\begin{thm}[\cite{EG:II}]
\label{thm_otcover}
Assume $\alpha$ is a positive contact form on a solid torus $V$ 
whose Reeb field is tangent to the boundary $\del V$. Choose 
any transverse curve $\gamma$ on the boundary torus $\del V$
which bounds a meridional disc in $V$. If the self-linking number
$s\ell k(\gamma)$ of this meridian is not equal to $-1$, then 
the pullback of $\alpha$ under some finite cover is an overtwisted 
contact form.
\end{thm}

It known from the inequality of \cite{Ben83} that if the initial 
contact structure is tight, the self linking number must 
satisfy $s\ell k\leq-\chi(D)=-1$. Thus, any self-linking number 
greater than $-1$ automatically implies an overtwisted structure
(which is of course preserved under covers). The nontrivial result
of this theorem is that for a tight structure, a self-linking number 
less than $-1$ implies an overtwisted cover.  
The techniques used in the proof of this theorem are a combination 
of perturbing characteristic foliations, manipulating singularities
of characteristic foliations, and using dynamical properties
of the characteristic foliation on $\del V$. 

\begin{remark}
It is necessary to distinguish between positive and negative contact
structures. In the case of a negative contact form (one for which 
$\alpha\wedge d\alpha<0$) on an invariant solid torus, the 
structure possesses an overtwisted cover if and only if the self-linking
number of a transverse meridian is not equal to $+1$. 
It is an easy exercise to show that the sign of the contact
form dual to a curl eigenfield is precisely the sign of the
eigenvalue. 
\end{remark}

>From these ingredients the following index may be defined:

\begin{dfn}
\label{def_Index}
Given a nonsingular Beltrami field $X$ on an invariant Riemannian 
solid torus $V$, define the index $\Index$ as follows.
\begin{enumerate}
\item
	If the eigenvalue $\lambda$ of $X$ with respect to the curl operator
	is zero, define $\Index := 0$. 
\item
	Otherwise, consider the characteristic foliation $(\del V)_\xi$
	of the contact structure orthogonal to $X$ on $\del V$. 
	If possible, choose $\gamma$ any meridional curve on $\del V$ 
	transverse to $\del V_\xi$ and define 
\begin{equation}
	\Index := \sign(\lambda)\left(s\ell k(\gamma)\right)+1 .
\end{equation}
\item
	If no transverse curve $\gamma$ exists, 
	define $\Index := \sign(\lambda)$.
\end{enumerate}		
\end{dfn}

\begin{thm}
\label{thm_Main}
Any $C^2$ or smoother nonvanishing Beltrami field on an invariant 
Riemannian solid torus $V$ having nonzero index $\Index$ possesses 
a contractible closed orbit. 
\end{thm}
{\em Proof:}
Assume that $\lambda>0$ as the negative case follows similarly. Since
$\Index\neq 0$, we are either in the case where there is 
a transverse meridian on the boundary of $V$ with self-linking 
number not equal to $-1$, or there is no transverse meridian. 
In the case where the transverse meridian exists, some finite
cover of the Beltrami field has a contractible closed orbit 
which is preserved by the covering projection. 

\begin{figure}[ht]
	{\epsfxsize=4.0in\epsfbox{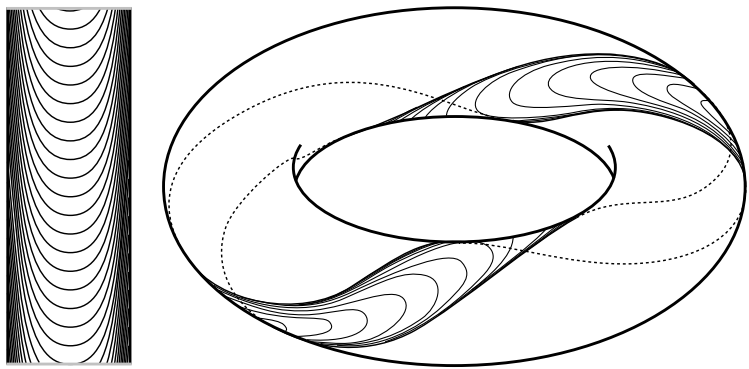}}
	\caption{A Reeb component in a two-dimensional foliation is a foliation 
	of an annulus by curves which limit onto the boundary components
	as illustrated [left, identify top and bottom]. 
	On the boundary torus $\del V$, if the 
	characteristic foliation $(\del V)_\xi$ possesses a Reeb 
	component [right], then there does not exist a transverse
	meridional curve.}
	\label{fig_reebcomp}
\end{figure}

If there does not exist a closed transversal, then a basic 
result in foliation theory implies that either 
(1) the characteristic foliation $(\del V)_\xi$ is entirely by 
meridional curves; or (2) the characteristic
foliation $(\del V)_\xi$ possesses a {\em Reeb component}, 
illustrated in Figure~\ref{fig_reebcomp}. In the former case, the 
contact structure is clearly overtwisted: any meridional disc 
in general position which 
spans one of these meridional curves has this boundary curve as a
limit cycle in its characteristic foliation. The existence of a 
periodic orbit then follows as earlier. 

In the latter case, where
a Reeb component exists in $(\del V)_\xi$, the following argument 
eliminates this possibility as a Beltrami field.  
Consider the Reeb field $Z$ associated to the contact form 
$\alpha$ dual to $X$. Since $Z$ preserves the contact structure 
$\xi$ and the boundary tours $\del V$, it must likewise preserve the
characteristic foliation $(\del V)_\xi$. Thus, if $(\del V)_\xi$ 
contains a closed curve, then, since $Z$ is everywhere transverse
to $(\del V)_\xi$, the entire boundary torus is swept out by forward
images of this curve under the flow of $Z$, and $(\del V)_\xi$ is 
a foliation by closed curves. Thus, a Reeb component (which always
possesses both closed and open curves as in the illustration) cannot 
arise as the characteristic foliation on an invariant solid torus. 
\qed

\section{Computation of $\Index$}
\label{sec_compute}

Given a contact structure $\xi$ on a 3-manifold $M$, the
determination of whether it is a tight or overtwisted
structure is a difficult question in general. By the 
Darboux Theorem (see, \eg, \cite{Aeb94,MS95}), every contact 
structure in dimension three is locally equivalent to 
the kernel of $dz+x\,dy$, 
which is a tight structure \cite{Ben83}. Thus, on the 
one hand, the property of being overtwisted is a 
decidedly global feature. However, the process of 
{\em Lutz twisting} \cite{Lut77} allows one to change a
tight structure into an overtwisted structure by means of
a $C^0$ alteration on an arbitrarily small open set in $M$. 
Thus, given a Beltrami field $X$ on $M$, determining 
whether the contact structure orthogonal to $X$ is overtwisted
is computationally intractable. Determining whether the 
universal cover is overtwisted is no less difficult.

However, to compute the index $\Index$, one does not need
information about the vector field on the entire three-dimensional
regime, but rather on some (arbitrary) two-dimensional 
meridional disc. We outline a method for easily computing 
$\Index$ from a $C^1$ approximation to $X$ along a two-dimensional
slice of the solid torus.

\begin{figure}[ht]
	{\epsfxsize=4.5in\epsfbox{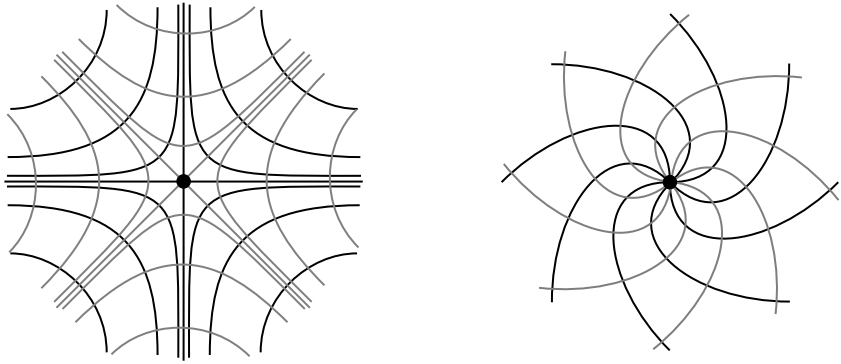}}
	\caption{Taking the orthogonal line field $X_D$ (grey) to the 
	characteristic foliation $D_\xi$ (black) leaves the index
	invariant.}
	\label{fig_duals}
\end{figure}

Choose a meridional disc $D$ with boundary curve 
$\gamma$. Orient $\gamma$ so that the contact form $\alpha:=g(X,\cdot)$ 
evaluates to a positive number on the tangents to $\gamma$: \ie, 
$\gamma$ points in roughly the same direction as $X$. This orientation 
on $\gamma$ induces in the usual way an orientation on the disc $D$. 
Denote by $X_D$ the projection of the vector field $X$ on to 
$D$ by orthogonal projection onto the tangent planes (orthogonal 
with respect to the metric $g$). 
\begin{prop}
If the vector field $X_D$ is generic (possesses a finite number
of nondegenerate rest points), the index $\Index$ of the Beltrami 
field $X$ can be computed by 
\begin{equation}
\label{eq_Index}
	\Index = 
\sign(\lambda)\left(1+\sum_{p : X_D(p)=0}\sigma(p)\Ind(X_D;p)\right)
\end{equation}
where for every rest point $p$ of $X_D$, $\sigma(p)$ is defined 
to be the sign $(+/-)$ of the dot product of $X(p)$ with the 
positive normal vector to $D$ at $p$, and the term $\Ind(X_D;p)$
denotes the standard Euler-Poincar\'e index of the planar vector field 
$X_D$ at $p$.
\end{prop}
{\em Proof:}
Assume that the characteristic foliation $D_\xi$ is generic in the 
above sense: there are a finite number of singular points $p$ at
which the contact structure $\xi$ is tangent to $D$, and the
characteristic foliation about these points appears locally as a 
source/sink or a saddle.
Following \cite{Ben83,Eli93}, there is a standard formula for 
computing the self-linking number of the transverse curve $\gamma=\del D$:
\begin{equation}
	s\ell k(\gamma) = 
	\sum_{p : D_\xi(p)=0}{\mbox{sign}}(T_pD,\xi_p)\Ind(D_\xi;p)
\end{equation}
Here, the sign of $(T_pD,\xi_p)$ is $+1$ when the orientations on 
the contact plane $\xi$ and the orientations on the tangent plane to 
disc $D$ at $p$ agree. Otherwise the sign is $-1$. However, we
only know the Beltrami field $X$ and the disc $D$ --- not the 
characteristic foliation. Determining the indices of the rest points 
of the characteristic foliation $D_\xi$ is accomplished via 
the projected field $X_D$ as follows. Since the contact structure
for $X$ is the plane field $\xi$ orthogonal to $X$, the characteristic
foliation at every point $q\in D$ is given by
\[
	D_\xi(q) := \xi_q\cap T_qD = (X_D(q))^\perp , 
\]
the line field orthogonal to $X_D$ at $q$. The proposition 
is proved by noting (1) fixed points of $X_D$ occur exactly
at fixed points of $D_\xi$; (2) the Euler-Poincar\'e 
index of $D_\xi$ at a fixed point $p$ is unchanged 
by looking at the orthogonal vector field,
as illustrated in Figure~\ref{fig_duals}. 
\qed

Computationally, this is extremely simple as $D$ can be chosen 
almost arbitrarily (one presumably chooses a 
$D$ which is ``nice'' in coordinates) and information about $X$
is required only on $D$ itself. 
The local index calculation is the most delicate portion of
the computation: the location of the orthogonal point $p$ 
is easy and the sign $\sigma(p)$ merely measures whether $X$ agrees
with the oriented normal to $D$, which is trivial to determine. 

\begin{figure}[ht]
	{\epsfxsize=3.5in\epsfbox{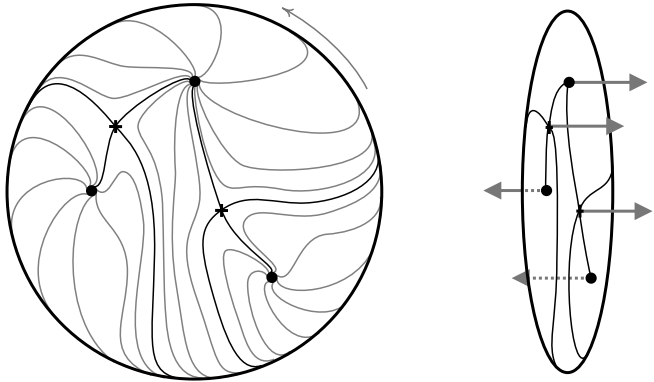}}
	\caption{An example of an $s\ell k$ computation given the 
	characteristic foliation on an oriented spanning disc $D$. 
	[left] $D_\xi$;	[right] the directions of the field $X$, positive 
	to the left and negative to the right. The self-linking 
	number is $s\ell k = -3$. Thus the Beltrami field $X$ has 
	index $\Index = -2$.}
	\label{fig_slkexample}
\end{figure}

An example of a characteristic foliation $D_\xi$ and the resulting
self-linking number is given in Figure~\ref{fig_slkexample}.

\section{Conclusions}
\label{sec_conc}

The index $\Index$ is an unusual object in that one inputs information 
about the vector field which is strictly {\em two}-dimensional, yet one
obtains data about the dynamics which is fully {\em three}-dimensional.
In Figure~\ref{fig_slkexample}, the {\em only} information about
the Beltrami field known is that (1) it is orthogonal to the disc $D$ at
five points; (2) locally near those five points the projected field 
is either source/sink or saddle type; and (3) at those 
five points the field points out or in as illustrated. From this 
data, it is inevitable that {\em somewhere} in the flow there exists 
a contractible closed flowline. 
This is a corollary of the contact-topological methods used --- the 
moral of the story is that the Beltrami condition hides
within it certain constraints on the dynamics which couple 
the dynamics of the vector field to the topology of the orthogonal
plane field. 

A deficiency of our theory is that it is not sharp. It is certainly 
possible for a Beltrami field to have the value $\Index=0$ and
yet still have a contractible periodic orbit: indeed, any tight contact
structure which has been Lutz twisted along a contractible closed curve
necessarily has trivial index as well as a contractible orbit. 
It appears certain (due to this mechanism of 
arbitrarily small Lutz twists) that no completely sharp 
computable index can be defined. What $\Index$ does, however, 
is detect if there are contractible orbits forced by the presence 
of non-localized overtwisted discs, and in this regime it is efficacious. 

It would be interesting to find a sharp lower bound on the number of 
periodic orbits present in the case of $\Index\neq 0$
(cf. the recent body of theory surrounding the 
{\em contact homology} of Eliashberg, Givental, and Hofer \cite{EGH00}).

\bibliographystyle{plain}


\end{document}